\documentclass[fleqn]{mat01}
\usepackage{times,mathtimy,amssymb,latexsym}
\begin{document}

\setcounter{page}{221} \firstpage{221}

\newtheorem{theoree}{Theorem}
\renewcommand\thetheoree{\Alph{theoree}}
\newtheorem{theor}[theoree]{\bf Theorem}
\newtheorem{theore}{\bf Theorem}
\newtheorem{lem}{Lemma}
\newtheorem{rema}{Remark}
\newtheorem{prob}{Problem}
\newtheorem{coro}{\rm COROLLARY}

\title{On $\pmb{p}$-quermassintegral differences function}

\markboth{Zhao Changjian and Wingsum Cheung}{On
$p$-quermassintegral differences function}

\author{ZHAO CHANGJIAN$^{1,*}$ and WINGSUM CHEUNG$^{2}$}

\address{$^{1}$Department of Information and Mathematics Sciences, College of
Science, China~Jiliang University, Hangzhou 310018, People's
Republic of China\\
\noindent $^{2}$Department of Mathematics, The University of Hong
Kong, Pokfulam Road, Hong~Kong\\
\noindent $^{*}$Author to whom correspondence should be
addressed.\\
\noindent E-mail: chjzhao@163.com; wscheung@hku.hk}

\volume{116}

\mon{May}

\parts{2}

\pubyear{2006}

\Date{MS received 3 August 2005; revised 25 December 2005}

\begin{abstract}
In this paper we establish Minkowski inequality and
Brunn--Minkowski inequality for $p$-quermassintegral differences
of convex bodies. Further, we give Minkowski inequality and
Brunn--Minkowski inequality for quermassintegral differences of
mixed projection bodies.
\end{abstract}


\maketitle

\vspace{6pt}
\section{Introduction}

The well-known classical Brunn--Minkowski inequality can be stated
as follows:

If $K$ and $L$ are convex bodies in $R^{n}$, then
\begin{equation}
V(K+L)^{1/n}\geq V(K)^{1/n}+V(L)^{1/n},
\end{equation}
with equality if and only if $K$ and $L$ are homothetic.

The Brunn--Minkowski inequality, has in recent decades,
dramatically extended its influence in many areas of mathematics.
Various applications have surfaced, for example, to probability
and multivariate statistics, shape of crystals, geometric
tomography, elliptic partial differential equations, and
combinatorics (see \cite{1,5,9,10,17}). Several remarkable analogs
have been established in other areas, such as potential theory and
algebraic geometry (see \cite{3,4,6,7,12,16}). Reverse forms and
similar forms of the inequality are important in the local theory
of Banach space (see \cite{17,18,20,21,22}). An elegant survey on
this inequality is provided by Gardner (see \cite{11}).

In fact, let $K$ and $L$ be convex bodies in $R^{n}$ and let
$0\leq i\leq n-1$. The Brunn--Minkowski inequality for
quermassintegral is the following inequality \cite{11}:
\begin{equation}
W_{i}(K+L)^{1/(n-i)}\geq W_{i}(K)^{1/(n-i)}+W_{i}(L)^{1/(n-i)}
\end{equation}
with equality if and only if $K$ and $L$ are homothetic.

Recently, $i$-{\it quermassintegral difference function} was
defined by Leng \cite{13} as follows:
\begin{equation*}
Dw_{i}(K,D)=W_{i}(K)-W_{i}(D), \quad \!K,D\in {\cal K}^{n},~
D\subset K ~{\rm and}~ 0\leq i\leq n-\!1.
\end{equation*}

Moreover, inequality (2) was extended to quermassintegral
differences of convex bodies  as follows \cite{13}:

\begin{theor}[\!]
If $K,L,$ and $D$ are convex bodies in $R^{n},$ $D\subset K$ and
$D'$ is a homothetic copy of $D,$ then
\begin{align}
Dw_{i}(K+L, D+D')^{1/(n-i)}\geq Dw_{i}(K,D)^{1/(n-i)}+
Dw_{i}(L,D')^{1/(n-i)},\hskip -1pc\phantom{00}
\end{align}
with equality for $0\leq i<n-1$ if and only if $K$ and $L$ are
homothetic and $(W_{i}(K),W_{i}(D))=\mu(W_{i}(L),W_{i}(D')),$
where $\mu$ is a constant.
\end{theor}

In \cite{8}, Firey introduced, for each real $p\geq 1$, new linear
combinations of convex bodies: For $K,L\in {\cal K}^{n}$, and
$\lambda, \mu\geq 0$ (both are not zero), {\it the Firey
combination}, $\lambda\cdot K+_{p}\mu\cdot L$, is a convex body.
The main aim of this paper is to establish the Brunn--Minkowski
inequality for quermassintegral differences about the Firey
combination, which is an extension of the inequality (3).

\begin{theore}[\!]
If $K,L,$ and $D$ are convex bodies in $R^{n},$ $D\subset K$ and
$D'$ is a homothetic copy of $D,$ then for $p\geq 1,$
\begin{align}
Dw_{i}(K\!+_{p}L, D\!+_{p}D')^{p/(n-i)}\geq
Dw_{i}(K,D)^{p/(n-i)}\!+\!Dw_{i}(L,D')^{p/(n-i)},\hskip
-1pc\phantom{00}
\end{align}
with equality for $0\leq i<n-p$ if and only if $K$ and $L$ are
homothetic $(p=1)$ {\rm (}or are dilates $(p>1))$ and
$(W_{i}(K),W_{i}(D))=\mu(W_{i}(L),W_{i}(D')),$ where $\mu$ is a
constant.
\end{theore}

For two convex bodies $K$ and $L$, an important inequality of
mixed volume is the well-known Minkowski inequality
\begin{equation*}
V_{1}(K,L)^{n}\geq V(K)^{n-1}V(L),
\end{equation*}
with equality if and only if $K$ and $L$ are homothetic.

In 1984, the inequality was extended to compact domains by Zhang
\cite{19} as\break follows:

\begin{theor}[\!]
If $K$ is a compact domain with piecewise $C^{1}$ boundary
$\partial K${\rm ,} and $L$ is a convex body in $R^{n}${\rm ,}
then
\begin{equation}
V_{1}(K,L)^{n}\geq V(K)^{n-1}V(L),
\end{equation}
with equality if and only if $K$ and $L$ are homothetic.
\end{theor}

Recently, inequality (5) was extended to volume differences by
Leng \cite{13} as follows:

\begin{theor}[\!]
Suppose that $K$ and $D$ are compact domains{\rm ,} $L$ is a
convex body{\rm ,} and $D\subset K, D'\subset L$ and $D'$ is a
homothetic copy of $D$. Then
\begin{equation}
(V_{1}(K,L)-V_{1}(D,D'))^{n}\geq Dv(K,D)^{n-1}Dv(L,D'),
\end{equation}
with equality if and only if $K$ and $L$ are homothetic and
$(V(K),V(D))=\mu(V(L),V(D')),$ where $\mu$ is a constant.
\end{theor}

In \cite{14}, Lutwak introduced the {\it mixed} $p$-{\it
quermassintegrals} $W_{p,0}(K,L),$ $W_{p,1}(K,L)$, $\dots,$
$W_{p,n-1}(K,L)$, for $K,L\in {\cal K}^{n}$, and real number
$p\geq 1$. The next aim of this paper is to establish the
Minkowski inequality for mixed $p$-quermassintegral differences of
convex bodies. A~new generalization of the classical Minkowski
inequality is presented as\break follows:

\begin{theore}[\!]
Let $K,L,$ and $D$ be convex bodies in $R^{n}$ and $D\subset K,~
D'\subset L$ and $D'$ be a homothetic copy of $D$. Then for $p\geq
1,$
\begin{equation}
(W_{p,i}(K,L)-W_{p,i}(D,D'))^{n-i}\geq
Dw_{i}(K,D)^{n-i-p}Dw_{i}(L,D)^{p},
\end{equation}
with equality for $0\leq i<n-p$ if and only if $K$ and $L$ are
homothetic $(p=1)$ {\rm (}or are dilates $(p>1))$ and
$(W_{i}(K),W_{i}(D))=\mu(W_{i}(L),W_{i}(D')),$ where $\mu$ is a
constant.
\end{theore}

On the other hand, we establish Minkowski inequality and
Brunn--Minkowski inequality for quermassintegral differences of
mixed projection bodies, respectively. which are extensions of
Lutwak's results \cite{15}.

\begin{theore}[\!]
Let $K,L,$ and $D$ be convex bodies in $R^{n},$ $D\subset K$ and
$D'$ a homothetic copy of $D$. Then for $0\leq j<n-2,$
\begin{align}
&Dw_{i}(\Pi_{j}(K+L),\Pi_{j}(D+D'))^{1/(n-i)(n-j-1)}\nonumber\\[.2pc]
&\quad\ \geq Dw_{i}(\Pi_{j} K,\Pi_{j}
D)^{1/(n-i)(n-j-1)}+Dw_{i}(\Pi_{j} L,\Pi_{j}
D')^{1/(n-i)(n-j-1)},\nonumber\\
\end{align}
with equality for $0\leq i<n-1$ if and only if $K$ and $L$ are
homothetic and $(W_{i}(K),W_{i}(D))=\mu(W_{i}(L),W_{i}(D')),$
where $\mu$ is a constant.
\end{theore}

\begin{theore}[\!]
Let $K,L,$ and $D$ be convex bodies in $R^{n},$ $D\subset K$, $D'$
be a homothetic copy of $D,$ and $0\leq j<n-1$. Then
\begin{align}
Dw_{i}(\Pi_{j}(K,L),\Pi_{j}(D,D'))^{n-1}\geq Dw_{i}(\Pi K, \Pi
D)^{n-j-1}Dw_{i}(\Pi L,\Pi D')^{j},\hskip -1pc\phantom{00}
\end{align}
with equality for $0\leq i<n-1$ if and only if $K$ and $L$ are
homothetic.
\end{theore}

The above interrelated notations, definitions and background
materials are given in \S2.

\section{Definitions and preliminaries}

The setting for this paper is the $n$-dimensional Euclidean space
${\Bbb R}^{n}(n>2)$. Let ${\cal K}^{n}$ denote the set of convex
bodies (compact, convex subsets with non-empty interiors) in
${\Bbb R}^{n}$.

Firey \cite{8} introduced, for each real $p\geq 1$, new linear
combinations of convex bodies: For $K,L\in {\cal K}^{n}$, and
$\lambda, \mu\geq 0$ (both are not zero), {\it the Firey
combination}, $\lambda\cdot K+_{p}\mu\cdot L$, is a convex body
whose support function is defined by
\begin{equation*}
h(\lambda\cdot K+_{p}\mu\cdot L,\cdot)^{p}=\lambda
h(K,\cdot)^{p}+\mu h(L,\cdot)^{p}.
\end{equation*}
Obviously, $\alpha\cdot K=\alpha^{1/p}K.$

The mixed quermassintegral
$W_{0}(K,L),W_{1}(K,L),\dots,W_{n-1}(K,L)$ of $K,L\in {\cal
K}^{n}$ are defined by
\begin{equation}
(n-i)W_{i}(K,L)=\lim_{\varepsilon\rightarrow
0}\frac{W_{i}(K+\varepsilon L)-W_{i}(K)}{\varepsilon},
\end{equation}
where
\begin{equation*}
W_{i}(K,L)=V(\underbrace{K,\dots,K}_{n-i-1},\underbrace{B,\dots,B}_{i},L).
\end{equation*}

The mixed $p$-quermassintegrals $W_{p,0}(K,L),
W_{p,1}(K,L),\dots,W_{p,n-1}(K,L)$, for $K,L\in {\cal K}^{n}$, and
real $p\geq 1$, are defined by \cite{14}
\begin{equation}
\frac{n-i}{p}W_{p,i}(K,L)=\lim_{\varepsilon\rightarrow
0}\frac{W_{i}(K+_{p}\varepsilon\cdot
L)-W_{i}(K)}{\varepsilon}.
\end{equation}

Of course for $p=1$, the mixed $p$-quermassintegral $W_{p,i}(K,L)$
is just $W_{i}(K,L)$. Obviously, $W_{p,i}(K,K)=W_{i}(K),$ for all
$p\geq 1$.

If $K_{1},\dots,K_{r}\in {\cal K}^{n}$ and
$\lambda_{1},\dots,\lambda_{r}\geq 0$, then the projection body of
the Minkowski linear combination
$\lambda_{1}K_{1}+\cdots+\lambda_{r}K_{r}\in {\cal K}^{n}$ can be
written as a symmetric homogeneous polynomial of degree $(n-1)$ in
$\lambda_{i}$ \cite{15}:
\begin{equation}
\Pi(\lambda_{1}K_{1}+\cdots+\lambda_{r}K_{r})=\sum \lambda_{i_{1}}
\dots\lambda_{i_{n-1}}\Pi_{i_{1}\cdots i_{n-1}},
\end{equation}
where the sum is a Minkowski sum taken over all $(n-1)$-tuples
$(i_{1},\dots,i_{n-1})$ of positive integers not exceeding $r$.
The body $\Pi_{i_1\dots i_{n-1}}$ depends only on the bodies
$K_{i_{1}},\dots,K_{i_{n-1}}$, and is uniquely determined by (12).
It is called {\it the mixed projection bodies} of
$K_{i_{1}},\dots,K_{i_{n-1}}$, and is written as
$\Pi(K_{i},\dots,K_{i_{n-1}}).$ If $K_{1}=\cdots=K_{n-1-i}=K$ and
$K_{n-i}=\cdots=K_{n-1}=L,$ then
$\Pi(K_{i_{1}},\dots,K_{i_{n-1}})$ will be written as
$\Pi_{i}(K,L).$ If $L=B$, then $\Pi_{i}(K,L)$ is denoted by
$\Pi_{i} K$ and when $i=0$, $\Pi_{i} K$ is denoted by $\Pi K.$

\section{Some lemmas}

The following results will be required to prove our main theorems.

\begin{lem}\hskip -.3pc {\rm \cite{14}.} \ \ If $p\geq 1,$ and $K,L\in
{\cal K}^{n},$ when $0\leq i<n,$ then
\begin{equation}
W_{i}(K+_{p}L)^{p/(n-i)}\geq W_{i}(K)^{p/(n-i)}+W_{i}(L)^{p/(n-i)},
\end{equation}
for $p>1$ with equality if and only if $K$ and $L$ are dilates{\rm
;} for $p=1$ with equality if and only if $K$ and $L$ are
homothetic.
\begin{equation}
W_{p,i}(K,L)^{n-i}\geq W_{i}(K)^{n-i-p}W_{i}(L)^{p}
\end{equation}
for $p>1$ with equality if and only if $K$ and $L$ are dilates{\rm
;} for $p=1$ with equality if and only if $K$ and $L$ are
homothetic.
\end{lem}

\begin{lem}\hskip -.3pc {\rm \cite{15}.} \ \ If $K,L\in {\cal K}^{n},$ and
$0\leq i<n,$ then
\begin{align}
&W_{i}(\Pi_{j}(K+L))^{1/(n-i)(n-j-i)}\nonumber\\[.2pc]
&\quad\ \geq W_{i}(\Pi_{j} K)^{1/(n-i) (n-j-1)} +W_{i}(\Pi_{j}
L)^{1/(n-i)(n-j-1)},\hskip -1pc\phantom{00}
\end{align}
with equality if and only if $K$ and $L$ are homothetic.\pagebreak
\begin{equation}
W_{i}(\Pi_{j}(K,L))^{n-1}\geq W_{i}(\Pi K)^{n-j-1}W_{i}(\Pi L)^{j},
\end{equation}
with equality if and only if $K$ and $L$ are homothetic.
\end{lem}

\begin{lem}\hskip -.3pc {\rm \cite{2} ({\it the Bellman's inequality}).} \ \ Let
$a=\{a_{1},\dots,a_{n}\}$ and $b=\{b_{1},\dots,b_{n}\}$ be two
series of positive real numbers and $p>1$ such that $a_{1}^{p}
-\sum_{i=2}^{n}a_{i}^{p}>0$ and $b_{1}^{p}
-\sum_{i=2}^{n}b_{i}^{p}>0$. Then
\begin{align}
\left(a_{1}^{p}-\!\sum_{i=2}^{n}a_{i}^{p}\right)^{1/p}\!+\!\left(
b_{1}^{p}-\!\sum_{i=2}^{n}b_{i}^{p}\right)^{1/p}
\leq\left((a_{1}+b_{1})^{p}-\!\sum_{i=2}^{n}(a_{i}+b_{i})^{p}\right)^{1/p}\hskip
-1pc\phantom{00}
\end{align}
with equality if and only if $a=\upsilon b$ where $\upsilon$ is a
constant.
\end{lem}

\begin{lem} If $a,b,c,d>0, 0<\alpha<1, 0<\beta<1$
and $\alpha+\beta=1$. Let $a>b$ and $c>d,$ then
\begin{equation}
a^{\alpha}c^{\beta}-b^{\alpha}d^{\beta}\geq (a-b)^{\alpha}(c-d)^{\beta},
\end{equation}
with equality if and only if $a/b=c/d.$
\end{lem}

\begin{proof}
Consider the following function
\begin{equation*}
f(x)=x^{\alpha}c^{\beta}-(x-b)^{\alpha}(c-d)^{\beta}, \quad x>0.
\end{equation*}
Let
\begin{equation*}
f'(x)=\alpha
c^{\beta}x^{\alpha-1}-\alpha(c-d)^{\beta}(x-b)^{\alpha-1}=0.
\end{equation*}
We get $x={bc}/{d}.$

On the other hand, if $x\in (0,\frac{bc}{d}),$ then $f'(x)<0$; if
$x\in (\frac{bc}{d},+\infty),$ then $f'(x)>0,$ and it follows that
\begin{equation*}
\min_{x>0}\{f(x)\}=f \left(\frac{bc}{d}
\right)=b^{\alpha}d^{\beta}.
\end{equation*}

This completes the proof.\hfill $\Box$
\end{proof}

\section{Inequalities for mixed $\pmb{p}$-quermassintegral
differences of convex bodies}

\setcounter{theore}{0}

\begin{theore}[\!] If $K,L,$ and $D$ are convex bodies in $R^{n},$
$D\subset K$ and $D'$  is a homothetic copy of $D,$ then for
$p\geq 1${\rm ,}
\begin{align}
Dw_{i}(K\!+_{p}L, D+_{p}D')^{p/(n-i)}\geq
Dw_{i}(K,D)^{p/(n-i)}+\!Dw_{i}(L,D')^{p/(n-i)},\hskip
-1pc\phantom{00}
\end{align}
with equality for $0\leq i<n-p$ if and only if $K$ and $L$ are
homothetic $(p=1)$ {\rm (}or are dilates $(p>1))$ and
$(W_{i}(K),W_{i}(D))=\mu(W_{i}(L),W_{i}(D')),$ where $\mu$ is a
constant.
\end{theore}

Using inequality (13) and in view of the Bellman's inequality, we
get the above Brunn--Minkowski inequality for quermassintegral
differences of the Firey combination.\pagebreak

Taking $p=1$ in inequality (19), the inequality (19) changes to
inequality (3).

In the following, we will prove the Minkowski inequality for mixed
$p$-quermassintegral differences of convex bodies.

\begin{theore}[\!] If $K, L$ and $D$ are convex bodies in $R^{n},$
$D\subset K, D'\subset L$ and $D'$ is a homothetic copy of $D,$
then for $p\geq 1,$
\begin{equation}
(W_{p,i}(K,L)-W_{p,i}(D,D'))^{n-i}\geq
Dw_{i}(K,D)^{n-i-p}Dw_{i}(L,D)^{p},
\end{equation}
with equality for $0\leq i<n-p$ if and only if  $K$ and $L$ are
homothetic $(p=1)$ {\rm (}or are dilates $(p>1))$ and
$(W_{i}(K),W_{i}(D))=\mu(W_{i}(L),W_{i}(D')),$ where $\mu$ is a
constant.
\end{theore}

\begin{proof}
From inequality (14), we have
\begin{equation}
W_{p,i}(K,L)^{n-i}\geq W_{i}(K)^{n-i-p}W_{i}(L)^{p},
\end{equation}
for $p>1$ with equality if and only if $K$ and $L$ are dilates;
for $p=1$ with equality if and only if $K$ and $L$ are homothetic
and
\begin{equation*}
W_{p,i}(D,D')^{n-i}=W_{i}(D)^{n-i-p}W_{i}(D')^{p}.
\end{equation*}
Hence,
\begin{align*}
W_{p,i}(K,L)-W_{p,i}(D,D') &\geq
W_{i}(K)^{(n-i-p)/(n-i)}W_{i}(L)^{p/(n-i)}\\[.2pc]
&\quad\ -W_{i}(D)^{(n-i-p)/(n-i)} W_{i}(D')^{p/(n-i)}.
\end{align*}
On the other hand, in view of inequality (18), we have
\begin{align*}
&W_{i}(K)^{(n-i-p)/(n-i)}W_{i}(L)^{p/(n-i)}-W_{i}(D)^{(n-i-p)/(n-i)}
W_{i}(D')^{p/(n-i)}\\[.2pc]
&\quad\
\geq(W_{i}(K)-W_{i}(D))^{(n-i-p)/(n-i)}(W_{i}(L)-W_{i}(D')^{p/(n-i)},\tag{21$'$}
\end{align*}
with equality if and only if
$W_{i}(K)/W_{i}(D)=W_{i}(L)/W_{i}(D').$

Thus,
\begin{align*}
W_{p,i}(K,L)-W_{p,i}(D,D')\geq
Dw_{i}(K,D)^{(n-i-p)/(n-i)}Dw_{i}(L,D)^{p/(n-i)}.
\end{align*}
Combining equality conditions of inequalities (21) and (21$'$), it
shows that the equality holds for $0\leq i<n-p$ if and only if $K$
and $L$ are homothetic ($p=1$) (or are dilates ($p>1$)) and
$(W_{i}(K),W_{i}(D))=\mu(W_{i}(L),W_{i}(D'))$, where $\mu$ is a
constant.

The proof is complete.\hfill $\Box$
\end{proof}

Taking $p=1$ in inequality (20), we get the following result.

\begin{coro}$\left.\right.$\vspace{.5pc}

\noindent Suppose that $K, L$ and $D$ are convex bodies{\rm ,} and
$D\subset K, D'\subset L, D'$ is a homothetic copy of $D$. Then
\begin{equation}
(W_{i}(K,L)-W_{i}(D,D'))^{n-i}\geq
Dw_{i}(K,D)^{n-i-1}Dw_{i}(L,D),
\end{equation}
with equality for $0\leq i<n-1$ if and only if $K$ and $L$ are
homothetic and $(W_{i}(K),W_{i}(D))=\mu(W_{i}(L),W_{i}(D')),$
where $\mu$ is a constant.
\end{coro}

Taking $i=0$ in (22), it changes to inequality (6).\pagebreak

\begin{rema}{\rm
Let $p=1,i=0$, and $D$, $D'$ be a single point in (20). Then (20)
reduces to the classical Minkowski inequality. Hence, (20) is a
generalization of the classical Minkowski inequality.}
\end{rema}

\section{Inequalities for quermassintegral differences of mixed
projection bodies}

In this section, we first establish the Brunn--Minkowski
inequality for quermassintegral differences of mixed projection
bodies as follows:

\begin{theore}[\!]
Let $K,L,$ and $D$ be convex bodies in $R^{n},$ $D\subset K$ and
$D'$ a homothetic copy of $D$. Then for $0\leq j<n-2,$
\begin{align}
&Dw_{i}(\Pi_{j}(K+L),\Pi_{j}(D+D'))^{1/(n-i)(n-j-1)}\nonumber\\[.2pc]
&\qquad\ \geq Dw_{i}(\Pi_{j} K,\Pi_{j}
D)^{1/(n-i)(n-j-1)}+Dw_{i}(\Pi_{j} L,\Pi_{j}
D')^{1/(n-i)(n-j-1)},\nonumber\\
\end{align}
with equality for $0\leq i<n-1$ if and only if $K$ and $L$ are
homothetic and $(W_{i}(K),W_{i}(D))=\mu(W_{i}(L),W_{i}(D')),$
where $\mu$ is a constant.
\end{theore}

\begin{proof}
Applying inequality (15), we have
\begin{align}
&W_{i}(\Pi_{j}(K+L))^{1/(n-i)(n-j-1)}\nonumber\\[.2pc]
&\quad\ \geq W_{i}(\Pi_{j} K)^{1/(n-i)(n-j-1)}+W_{i}(\Pi_{j}
L)^{1/(n-i)(n-j-1)},\hskip -1pc\phantom{00}
\end{align}
with equality if and only if $K$ and $L$ are homothetic.
\begin{align}
&W_{i}(\Pi_{j} (D+D'))^{1/(n-i)(n-j-1)}\nonumber\\[.2pc]
&\quad\ =W_{i}(\Pi_{j} D)^{1/(n-i)(n-j-1)}+W_{i}(\Pi_{j}
D')^{1/(n-i)(n-j-1)}.
\end{align}

From (24) and (25), we obtain that
\begin{align}
&Dw_{i}(\Pi_{j}(K+L),\Pi_{j}(D+D'))\nonumber\\[.2pc]
&\quad\ \geq [W_{i}(\Pi_{j}
K)^{1/(n-i)(n-j-1)}+W_{i}(\Pi_{j} L)^{1/(n-i)(n-j-1)}]^{(n-i)(n-j-1)}\nonumber\\[.2pc]
&\qquad\ -[W_{i}(\Pi_{j}D)^{1/(n-i)(n-j-1)} +W_{i}(\Pi_{j}
D')^{1/(n-i)(n-j-1)}]^{(n-i)(n-j-1)}.
\end{align}

From (26) and in view of the Bellman's inequality,
\begin{align*}
&Dw_{i}(\Pi_{j}(K+L),\Pi_{j}(D+D'))^{1/(n-i)(n-j-1)}\\[.2pc]
&\quad\ \geq (W_{i}(\Pi_{j} K)-W_{i}(\Pi_{j}
D))^{1/(n-i)(n-j-1)}\\[.2pc]
&\qquad\ +(W_{i}(\Pi_{j} L)-W_{i}(\Pi_{j} D'))^{1/(n-i)(n-j-1)}.
\end{align*}

This completes the proof.\hfill $\Box$
\end{proof}

Taking $i=0,j=0$ in inequality (23), we obtain the following
result.

\begin{coro}$\left.\right.$\vspace{.5pc}

\noindent Let $K,L,$ and $D$ be convex bodies in $R^{n},$
$D\subset K$ and $D'$ is a homothetic copy of $D$. Then
\begin{align}
&Dv(\Pi(K+L),\Pi(D+D'))^{1/n(n-1)}\nonumber\\[.2pc]
&\quad\ \geq Dv(\Pi K,\Pi D)^{1/n(n-1)}+Dv(\Pi L,\Pi
D')^{1/n(n-1)}
\end{align}
with equality if and only if $K$ and $L$ are homothetic and
$V((K),V(D))=\mu(V(L),V(D')),$ where $\mu$ is a constant.
\end{coro}

This is just a projection form of `Theorem~1' which was given by
Leng \cite{13}.

\begin{rema}{\rm Let $D$ and $D'$ be a single point in (23). Then (23)
changes to (16). This shows that (23) is a generalization of the
Brunn--Minkowski inequality for mixed projection bodies.

In the following, we establish the Minkowski inequality for
quermassintegral differences of mixed projection bodies.}
\end{rema}

\begin{theore}[\!]
Let $K,L,$ and $D$ be convex bodies in $R^{n},$ $D\subset K,$ $D'$
is a homothetic copy of $D,$ and $0\leq j<n-1$. Then
\begin{align}
Dw_{i}(\Pi_{j}(K,L),\Pi_{j}(D,D'))^{n-1}\geq Dw_{i}(\Pi K, \Pi
D)^{n-j-1}Dw_{i}(\Pi L,\Pi D')^{j},\hskip -1pc\phantom{00}
\end{align}
with equality for $0\leq i<n-1$ if and only if $K$ and $L$ are
homothetic and $(W_{i}(K),W_{i}(D))=\mu(W_{i}( L),W_{i}(D'))$
where $\mu$ is a constant.
\end{theore}

\begin{proof}
Applying inequality (16), we have
\begin{equation*}
W_{i}(\Pi_{j}(K,L))\geq W_{i}(\Pi K)^{(n-j-1)/(n-1)}W_{i}(\Pi
L)^{j/(n-1)},
\end{equation*}
with equality if and only if $K$ and $L$ are homothetic.
\begin{equation*}
W_{i}(\Pi_{j}(D,D'))= W_{i}(\Pi D)^{(n-j-1)/(n-1)}W_{i}(\Pi
D')^{j/(n-1)}.
\end{equation*}
Hence, from inequality (18) in Lemma~4, we obtain that
\begin{align*}
&Dw_{i}(\Pi_{j}(K,L),\Pi_{j}(D,D'))\\[.2pc]
&\quad\ \geq W_{i}(\Pi K)^{(n-j-1)/(n-1)}W_{i}(\Pi
L)^{j/(n-1)}\\[.2pc]
&\qquad\ -W_{i}(\Pi D)^{(n-j-1)/(n-1)}W_{i}(\Pi
D')^{j/(n-1)}\\[.2pc]
&\quad\ \geq (W_{i}(\Pi K)-W_{i}(\Pi D))^{(n-j-1)/(n-1)}(W_{i}(\Pi
L)-W_{i}(\Pi D'))^{j/(n-1)}.
\end{align*}

The proof is complete.\hfill $\Box$
\end{proof}

Taking $i=0,j=1$ in inequality (28), inequality (28) changes to
the following result.

\begin{coro}$\left.\right.$\vspace{.5pc}

\noindent Let $K,L,$ and $D$ be convex bodies in $R^{n},$
$D\subset K$ and $D'$ a homothetic copy of $D$. Then
\begin{equation}
Dv(\Pi_{1}(K,L),\Pi_{1}(D,D'))^{n-1}\geq Dv(\Pi K,\Pi
D)^{n-2}Dv(\Pi L,\Pi D'),
\end{equation}
with equality for $0\leq i<n-1$ if and only if $K$ and $L$ are
homothetic and $(V(K),V(D))=\mu(V( L),V(D'))$ where $\mu$ is a
constant.
\end{coro}

This is just a projection form of inequality (6).

\begin{rema}{\rm Let $D$, $D'$ be a single point in (29). Then (29)
changes to the following inequality:

If $K,L\in {\cal K}^{n}$, then
\begin{equation*}
V(\Pi_{1}(K,L))^{n-1}\geq V(\Pi K)^{n-2}V(\Pi L)
\end{equation*}
with equality if and only if $K$ and $L$ are homothetic.

This is just {\it the Minkowski inequality for mixed projection
bodies} which was given by Lutwak \cite{15}.}
\end{rema}

\section{Two open problems}

In the following, we pose two open problems:

\begin{prob}{\rm
Let $K_{i}, i=1,2,\dots,n$ and $D_{i}, i=1,2,\dots,n$ be convex
bodies in $R^{n}$, $D_{i}\subset K_{i}$ and $D_{i}'$ a homothetic
copy of $D_{i}, i=1,2,\dots,n$, respectively. Then for $0\leq
r\leq n$,
\begin{align}
&(V(K_{1},\dots,K_{n})-V(D_{1},\dots,D_{n}))^{r}\nonumber\\[.2pc]
&\qquad
\geq\prod_{j=1}^{r}(V(\underbrace{K_{j},\dots,K_{j}}_{r},K_{r+1},
\dots,K_{n})
-\!V(\underbrace{D_{j},\dots,D_{j}}_{r},D_{r+1},\dots,D_{n})).\nonumber\\
\end{align}
}
\end{prob}

\begin{rema}{\rm In (30), taking $r=n$, we obtain that
\begin{equation}
(V(K_{1},\dots,K_{n})-V(D_{1},\dots,D_{n}))^{n}\geq\prod_{j=1}^{n}(V(K_{j})-V(D_{j})).
\end{equation}
Taking $K_{1}=\cdots=K_{n-1}=K, K_{n}=L, D_{1}=\cdots=D_{n-1}=D,
D_{n}=D'$ in (31), inequality (31) changes to
\begin{equation*}
(V_{1}(K,L)-V_{1}(D,D'))^{n}\geq Dv(K,D)^{n-1}Dv(L,D').
\end{equation*}
This is just inequality (6).

On the other hand, let $D$ and $D'$ be a single point in (30).
Then (30) changes to the well-known Aleksandrov--Fenchel
inequality.}
\end{rema}

\begin{prob}{\rm Let $K_{i}, i=1,\dots,n-1$ and
$D_{i},i=1,\dots,n-1$ be convex bodies in $R^{n}$, $D_{i}\subset
K_{i}$ and $D_{i}', i=1,\dots,n-1$ a homothetic copy of $D_{i}$,
respectively. Then for $0\leq r\leq n-1$,
\begin{align}
&Dv(\Pi(K_{1},\dots,K_{n-1}),\Pi(D_{1},\dots,D_{n-1}))^{r}\nonumber\\[.2pc]
&\ \
\geq\prod_{j=1}^{r}\!Dv(\Pi(\underbrace{K_{j},\dots,\!K_{j}}_{r},K_{r+1},\dots,\!K_{n}),
\Pi(\underbrace{D_{j},\dots,\!D_{j}}_{r},D_{r+1},\dots,\!D_{n})).
\end{align}}
\end{prob}

\begin{rema}{\rm In (32), taking $r=n-1$, we obtain that
\begin{equation}
Dv(\Pi(K_{1},\dots,K_{n-1}),\Pi(D_{1},\dots,D_{n-1}))^{n-1}\geq\prod_{j=1}^{n-1}Dv(\Pi
K_{j},\Pi D_{j}).
\end{equation}
Taking $K_{1}=\cdots=K_{n-2}=K, K_{n-1}=L, D_{1}=\cdots=D_{n-2}=D,
D_{n-1}=D'$ in (33), inequality (33) changes to
\begin{equation}
Dv(\Pi_{1}(K,L),\Pi_{1}(D,D'))^{n}\geq Dv(\Pi (K,D))^{n-1}Dv(\Pi (L,D')).
\end{equation}
This is just inequality (29) which is proved in this paper.

On the other hand, let $D$ and $D'$ be a single point in (32).
Then (32) changes to the well-known Aleksandrov--Fenchel
inequality for mixed projection bodies which was given by Lutwak
\cite{15}.}
\end{rema}

\section*{Acknowledgements}

The first author (ZC) was supported by National Natural Sciences
Foundation of China (10271071), Zhejiang Provincial Natural
Science Foundation of China (Y605065), Foundation of the Education
Department of Zhejiang Province of China (20050392) and the
Academic Mainstay of Middle-age and Youth Foundation of Shandong
Province of China (200203). The second author (WC) was partially
supported by the Research Grants Council of the Hong Kong SAR,
China (Project No: HKU7040/03P).

\end{document}